\documentclass[11pt, reqno]{article}
\usepackage{amsmath,amstext,amsthm,amsfonts}
\usepackage[dvips]{graphicx}
\usepackage{epic,eepic}

\theoremstyle{plain}
\newtheorem{theorem}{Theorem }[section]

\newtheorem{lemma}[theorem]{Lemma}
\newtheorem{corollary}[theorem]{Corollary}
\newtheorem{maintheorem}{Theorem}

\theoremstyle{definition}
\theoremstyle{remark}

\newtheorem{definition}[theorem]{Definition}

\newcommand{\field}[1]{\mathbb{#1}}
\newcommand{\real}{\field{R}}

\newcommand{\Circ}{\field{S}}

\newcommand{\al} {\alpha}

\newcommand{\de} {\delta}       \newcommand{\De}{\Delta}

\begin{document}

\title{A remark on capillary surfaces in a 3-dimensional space of constant
curvature}
\author{ Alexander Arbieto, Carlos Matheus and Marcos Petr\'ucio 
\thanks{The first and second authors were partially supported by Faperj-Brazil. The third author
was partially supported by CNPq-Brazil} }

\date{June 16, 2003}
\maketitle
\begin{abstract}
{ We generalize a theorem by J. Choe on capillary surfaces for arbitrary
3-dimensional spaces of constant curvature. The main tools in this paper are an
extension of a theorem of H. Hopf due to S.-S. Chern and two index lemmas
by J. Choe. }
\end{abstract}

\section{History}

	A well known theorem due H. Hopf \cite{H} state that a CMC immersion of the 
sphere into the Euclidean three-dimensional space is a round sphere.
In 1982, in Rio de Janeiro, in occasion of a International Congress at IMPA, S. S. Chern
\cite{C} showed a generalization of this theorem when the ambient space has constant
sectional curvature. Recently, J. Choe \cite{Ch}, using sufficient hypothesis,
generalized Hopf's Theorem for immersion of the closed disk in $R^3$. Following these
lines of ideias, the first result in this paper can be stated as:

\begin{maintheorem}\label{t.1} \textit{Let $S$ be a CMC immersed compact $C^{2+\al}$ 
surface of disk type in a 3-dimensional ambient space $M^3$ of constant
curvature ($C^{2+\al}$ surface means $C^{2+\al}$ up to and including
$\partial S$ and $\partial S$ is $C^{2+\al}$ up to and including its vertices).
Suppose that the regular components of $\partial S$ are lines of curvature. If
the number of vertices with angle $<\pi$ is less than or equal to $3$, then $S$ is totally
umbilic.}
\end{maintheorem}

This kind of theorem is motivated by the study of capillary surfaces. In fact,
J. Nitsche, in 1995, showed that a regular capillary immersion
(see definition~\ref{d.capillary}) of the closed disk in the sphere is either a plane
disk or a piece of a round sphere. In \cite{Ch} this result was obtained for capillary 
immersion without strong regularity assumptions. In 1997 Ros-Souam
\cite{RS} showed a version of Nitsche's theorem for ambient space with constant sectional curvature. They used Chern's
extension of Hopf's theorem. This motivated us to formulate the theorem below:

\begin{maintheorem}\label{t.2} \textit{Let $U\subset M^3$ be a domain of a $3$-dimensional
space of constant curvature bounded by totally umbilic surfaces. If $S$ is a
capillary surface in $U$ of disk type which is $C^{2+\al}$ and $S$ has less than
$4$ vertices with angle $<\pi$, then $S$ is totally umbilic. }
\end{maintheorem}

The paper is organized as follows. In section 2 we recall the context of Chern's
work~\cite{C}, state the main theorems needed here and briefly define the concept
of rotation index of the lines of curvature at umbilic points. Finally, we prove the
theorems A and B in the section 3.

\section{Some Lemmas}

In this section we fix some notation and briefly sketch the proof of the main
tools used here: S.-S. Chern's generalization of Hopf's theorem and the two
index lemmas by J. Choe.

Let $M^3$ be a 3-dimensional manifold of constant curvature $c$. Following
Chern~\cite{C}, if $X:S\rightarrow M$
is an immersed surface and $p\in S$, we can fix an orthornormal local frame
$e_1, e_2, e_3$ such that $e_3$ is the unit normal vector to $S$ at $x$,
supposing $S$ orientable ($x\in S$ is a point near to $p$). If $\theta_i$
denotes the coframe ($i=1,2,3$), then $\theta_3=0$. The first and second
fundamental forms are $I=\theta_1^2+\theta_2^2$ and $II=h_{11}\theta_1^2 +
2h_{12}\theta_1\theta_2 +h_{22}\theta_2^2$, respectively. 

Recall that the invariants $H=\frac{1}{2}(h_{11}+h_{22})$ and $\widetilde{K}=h_{11}h_{22}-h_{12}^2$ are
the mean curvature and the total curvature of S, where $S$ has the induced
Riemannian metric. By the structure
equations (see~\cite{C}), we have that the Gaussian curvature is
$K=\widetilde{K} + c$. With this setting, we recall the definitions:

\begin{definition} $S$ is \emph{totally umbilical} (resp. \emph{totally
geodesic}) if $II-H\cdot I=0$ (resp. if $II=0$).
\end{definition} 

Defining $\phi =\theta_1 +i\theta_2$, we have a complex structure on $S$. Note
that $II-H\cdot
I=\frac{1}{2}(h_{11}-h_{22})(\theta_1^2-\theta_2^2)+2h_{12}\theta_1\theta_2$ is
a trace zero form and it is the real part of the complex $2$-form
$\Phi=\widetilde{H}\phi^2$, where
$\widetilde{H}=\frac{1}{2}(h_{11}-h_{22})-h_{12}i$. Since $\Phi$ is uniquely
determined by $II-H\cdot I$ and $II- H\cdot I$ is associated to
$S$, $\Phi$ is a globally defined $2$-form, independent of a choice of local
frames.

In the work~\cite{H}, H. Hopf shows that, in the case $M=\real^3$, i.e., $c=0$, if the mean curvature is constant, then
$\Phi$ is holomorphic on $S$. However, a more general fact is true, as proved by Chern:

\begin{lemma}[Theorem 1 of Chern~\cite{C}]\label{l.Chern} If $H\equiv const.$, then $\Phi$ is a
holomorphic $2$-form on $S$.
\end{lemma} 

The holomorphicity of $\Phi$ was used by Hopf to prove that an immersed sphere
$f:S^2\rightarrow\real^3$ of \emph{constant mean curvature} (CMC) is round. Indeed, this
follows from a standard result about Riemman surfaces which says that, except by
the trivial $2$-form $\Phi=0$, there is
no holomorphic $2$-form on a compact Riemman surface of zero genus. With the
same argument, as a corollary of~\ref{l.Chern}, Chern was able to conclude that:

\begin{corollary}[Theorem 2 of Chern~\cite{C}] If $f:S^2\rightarrow M^3$ is an
CMC immersed sphere and $M^3$ has constant mean curvature then $f$ is
totally
umbilic.  
\end{corollary}

On the other hand, in the case of surfaces with boundary (with ambient space
$M=\real^3$), Choe extends Hopf's arguments to study \emph{capillary surfaces}.
In order to make the ideias of Hopf works in his case, Choe introduce a natural
concept of \emph{rotation index} of the lines of curvature at umbilic points
(\emph{including boundary points}). Now, as a preliminary work, we consider
Choe's notion of rotation index in the context of a general ambient space $M^3$ of
constant curvature. 

Consider a point $p\in\partial S$. Let $\psi:D_h\rightarrow S$ be a conformal
parametrization of a neighborhood of $p$ in $S$, where $D_h=\{ (x,y)\in D :
y\geq 0\}$ is a half disk and the diameter $l$ of $D_h$ is mapped into $\partial
S$. Let $F$ be the line field on $D_h$ obtained by pulling back (by
$\psi$)
the lines of curvature of $S$. If $\psi(l)$ is a line of curvature of $S$, we can
extend $F$ to a line field on $D$ by reflection through the diameter $l$. For
simplicity, the extension of $F$ to $D$ is denoted by $F$. At this point, it is
natural define the rotation index of the lines of curvature at an umbilic point
$p\in\partial S$ to be \emph{half} of the index of $F$ at $\psi^{-1}(p)$. Clearly
this definition is independent of the choice of the parametrization $\psi$. However,
the definition only makes sense if we show that the umbilic points on 
$\partial S$ are isolated. But this fact follows from an easy argument:

Following Choe~\cite{Ch}, the equation of the lines of curvature, in complex
coordinates is given by 
$$\Im( \Phi ) =0 ,$$    
where $\Im z$ denotes the imaginary part of $z$.

So the rotation index of the lines of curvature is 
$$r=\frac{1}{2\pi}\de (\text{arg} \phi)=-\frac{1}{4\pi}\de
(\text{arg}\Phi),$$
where $\de$ is the variation as one winds once around an isolated umbilic point
$p$. In particular, if $p$ is an interior point of $S$ (i.e., $p\notin S$) and
is a zero of order $n$ of $\Phi$, then $\de(\text{arg}\Phi )=2\pi n$.
Consequently, 
\begin{equation}\label{pto interior}
r=-\frac{n}{2}\leq -\frac{1}{2}.
\end{equation}

Suppose now that $\Phi$ has a zero (resp., pole) of order $n>0$ (resp., $-n>0$)
at a boundary umbilic point $p$. Then, 
\begin{equation}\label{pto bordo} r=\frac{1}{2}\big[ -\frac{1}{4\pi}\de
(\text{arg}\Phi)\big] =-\frac{n}{4} 
\end{equation} 

With these equations in mind, Choe proves the following lemma, which compares
interior umbilic points and boundary umbilic points. 

\begin{lemma}\label{l.Choe} Let $S$ be a CMC immersed $C^{2+\al}$ surface (up to
and including the boundary $\partial S$). Suppose that $\partial S$ consist of
$C^{2+\al}$ curves (up to and including some possible singular points called
vertices). If the regular components of $\partial S$ are lines of curvature,
then:
\begin{enumerate}\item The boundary umbilic points of $S$ are isolated;

 		 \item The boundary umbilic points which are not vertices have,
		 at most, rotation index $-1/4$ ;
		 
		 \item The vertices of $S$ with angle $<\pi$ have rotation index
		 $\leq 1/4$ and the vertices with angle $>\pi$ have rotation
		 index $\leq -1/4$.
\end{enumerate}
\end{lemma}

The proof of this lemma is a straightforward consequence (with only minor
modifications) of Choe's proof of
lemma 2 in~\cite{Ch}. 

Now we are in position to prove the man results of this paper.

\section{Proof of the theorems}

%We are now in position to prove the first main result of this paper:

%\begin{maintheorem}\label{t.1} Let $S$ be a CMC immersed compact $C^{2+\al}$ 
%surface of disk type in a 3-dimensional ambient space $M^3$ of constant
%curvature ($C^{2+\al}$ surface means $C^{2+\al}$ up to and including
%$\partial S$ and $\partial S$ is $C^{2+\al}$ up to and including its vertices).
%Suppose that the regular components of $\partial S$ are lines of curvature. If
%the number of vertices with angle $<\pi$ is less than or equal to $3$, then $S$ is totally
%umbilic.
%\end{maintheorem}

\begin{proof}[Proof of theorem~\ref{t.1}] Fix $\psi:D\rightarrow S$ a conformal parametrization and $F$ the
pull-back under $\psi$ of the lines of curvature of $S$. Since $\partial S$ are
lines of curvature, we can apply the Poincar\'e-Hopf theorem (even in the
case that $\psi$ is a parametrization of a boundary point) to conclude
that, if the
number of singularities of $F$ is finite, the sum of rotation indices is equal
to 1. Let $A$ be the set of such singularities. Suppose that $A$ is finite.
Using equation~\ref{pto interior}, lemma~\ref{l.Choe} and, by hypothesis, the
number of vertices with angle is $<\pi$ is $\leq 3$, we get the estimate:
$$\sum r(p)\leq 3/4 ,$$ a contradiction with Poincar\'e-Hopf's theorem.

Therefore, $A$ is infinite. In particular, since the number of vertices is
finite, we have an infinite set $A-\{ \text{vertices} \}\subset \{ \text{zeros
of }\Phi \}$. But $\Phi$ is holomorphic. In particular, this implies that $A=S$,
so $S$ is totally umbilic.        
\end{proof}

We point out that, as remarked by Choe~\cite{Ch}, Remark 1, the condition on the
number of vertices with angle $<\pi$ is necessary. In fact, a rectangular region
in a cylinder $N=\Circ\times\real^1\subset\real^3$ bounded by two straight lines
and two circles provides a counter-example with 4 vertices with angle $\pi /2$
and rotation index $1/4$.

Before starting the proof of the second main result, we recall the definition:

\begin{definition}\label{d.capillary} A \emph{capillary surface} $S$ in a domain $U$ of a
$3$-dimensional space $M^3$ of constant curvature is a CMC immersed surface
which meets $\partial U$ along $\partial S$ at a constant angle.  
\end{definition}

As a immediate consequence of theorem~\ref{t.1}, we have the theorem~\ref{t.2}:

%\begin{maintheorem}\label{t.2} Let $U\subset M^3$ be a domain of a $3$-dimensional
%space of constant curvature bounded by totally umbilic surfaces. If $S$ is a
%capillary surface in $U$ of disk type which is $C^{2+\al}$ and $S$ has less than
%$4$ vertices with angle $<\pi$, then $S$ is totally umbilic. 
%\end{maintheorem}

\begin{proof}[Proof of theorem~\ref{t.2}] This follows from theorem~\ref{t.1}
and the Terquem - Joachimsthal
theorem~\cite{S} that says: 

 {\it ``If $C=S_1\cap S_2$ is a line of curvature of $S_1$, then
$C$ is also a line of curvature of $S_2$ if and only if $S_1$ intersect $S_2$ at
a constant angle along $C$.''}
\end{proof}

A result for capillary \emph{hypersurfaces} with the same flavor of theorem~\ref{t.2} 
was also obtained by Choe
(see~\cite{Ch}, theorem 3). However, these arguments does not work \emph{a
priori} for more general ambient spaces than $\real^n$ since the following
fact (valid only in $\real^n$) is used: if $X$ denotes the position vector on
$S$ from
a fixed point and \textbf{H} is the mean curvature vector, then $\De
X=\textbf{H}$. In particular, it is an open question if there exists
\emph{unbalanced} capillary hypersurfaces in the conditions of theorem 3 of
Choe.\footnote{J. Choe pointed out to one of the authors that, in fact, there exists a
generalization of theorem 3 of~\cite{Ch} to be published elsewhere.}

To finish this paper we point out that another kind of generalization of the theorem 3 of Choe cited above is obtained
by replacing
the mean curvature by the \emph{higher order curvatures $H_r$}. In this direction,
Choe showed in~\cite{Ch}, theorem $4$, that if an immersed hypersurface
$S\subset\real^{n+1}$ has constant mean
curvature and $H_r$ is constant for some $r\geq 2$, then $S$ is a hypersphere.
Moreover, we can replace the constancy of the mean curvature by $S$ is embedded,
as Ros proved~\cite{R}. Furthermore, for general ambients of constant curvature
and supposing only that $H_r/H_l$ is constant, Koh-Lee~\cite{KL} were able to
get the same result. Recently, Alencar-Rosenberg-Santos~\cite{ARS} proved a result in this direction with
an extra hypothesis on the Gauss image of $S$ (with ambient space
$\Circ^{n+1}$). 

\textbf{Acknowledgements.} The authors are grateful to professor Manfredo do Carmo for his
encouragement to write this short note. Finally, the authors are thankful to IMPA and his staff.

\bibliographystyle{alpha}
\bibliography{bib}

\noindent 
		\textbf{Alexander Arbieto} ( alexande{\@@}impa.br )\\
		\textbf{Carlos Matheus} ( matheus{\@@}impa.br )\\
		\textbf{Marcos Petr\'ucio} ( petrucio{\@@}impa.br )\\
		IMPA, Est. D. Castorina, 110, Jardim Bot\^anico, 22460-320\\
		Rio de Janeiro, RJ, Brazil

\end{document}